\documentclass[onefignum,onetabnum]{siamonline190516}

\usepackage[utf8]{inputenc}
\usepackage[english]{babel}
\usepackage[titletoc]{appendix}
\usepackage{dsfont}
\usepackage{amsfonts}
\usepackage{bbm}
\usepackage{graphicx}
\usepackage{pgf,tikz}
\usepackage{nicefrac}
\usepackage{cleveref}
\usepackage{enumitem}
\setlist[enumerate]{leftmargin=.5in}
\usepackage{algorithmic}
\usepackage{todonotes}
\usepackage{placeins}

\numberwithin{equation}{section}
\definecolor{mygreen}{RGB}{22, 204, 70}

\newsiamthm{assumption}{Assumption}
\newsiamremark{remark}{Remark}
    \renewenvironment{equation*}{\[}{\]\ignorespacesafterend}
\def\Uad{{\mathcal P}}
\def\Vad{{\mathcal{X}}}
\def\N{\mathbb{N}}
\def \R{\mathbb{R}}
\def\eps{\varepsilon}
\def\P{\mathcal{P}}
\def\supp{\textnormal{supp}}
\def\Prob{\mathbb P}
\def\Ew{\mathbb E}
\def\mX{\mathcal X}
\def\mY{\mathcal Y}

\def\X{X}

\def\x{x}

\def\din{{d_{\textnormal{par}}}}

\def\dpar{{d_{\textnormal{des}}}}
\def\U{\Theta}
\def\u{\theta}

\def\Z{Z}
\def\mP{\Uad} % Parameter space
\DeclareMathOperator*{\argmin}{arg\,min}

\newcommand{\Proj}{{\mathrm{Proj}_{\Uad}}}

\newcommand{\cu}{{\underline{\mathsf S}}}
\newcommand{\co}{{\overline{\mathsf S}}}
\newcommand{\nX}{{_{\mX}}}
\newcommand{\nY}{{_{\mY}}}
\newcommand{\nT}{{_{\Uad}}}
\renewcommand{\Delta}{{\boldsymbol{\delta}}}
\newcommand{\CD}{\mathsf C_{_\Delta}}
\newcommand{\LD}{\mathsf L_{_\Delta}}

\newcommand{\LZ}{\mathsf L_{_Z}}
\newcommand{\Me}{\mathsf M_{_1}}
\newcommand{\Mz}{\mathsf M_{_2}}
\newcommand{\D}{\mathsf D}

\author{Lukas Pflug\thanks{Central Institute for Scientific Computing (ZISC), \email{lukas.pflug@fau.de}}%\thanks{Department of Mathematics, Chair of Applied Mathematics (Continuous Optimization), Friedrich-Alexander University Erlangen-Nürnberg (FAU),}
    \and
    Max Grieshammer\thanks{Department of Mathematics, Chair of Applied Mathematics (Continuous Optimization), Friedrich-Alexander University Erlangen-Nürnberg (FAU), \email{max.grieshammer@fau.de}, \email{andrian.uihlein@fau.de}, \email{michael.stingl@fau.de}}
    \and
    Andrian Uihlein\footnotemark[2]
    \and
    Michael Stingl\footnotemark[2]
    }

\title{CSG: A stochastic gradient method for a wide class of optimization problems appearing in a machine learning or data-driven context}

\begin{document}

\maketitle

\begin{abstract}
    A recent article introduced the \emph{continuous stochastic gradient method} (CSG) for the efficient solution of a class of stochastic optimization problems. While the applicability of known stochastic gradient type methods is typically limited to expected risk functions, no such limitation exists for CSG. This advantage stems from the computation of design dependent integration weights, allowing for optimal usage of available information and therefore stronger convergence properties. However, the nature of the formula used for these integration weights essentially limited the practical applicability of this method to problems in which stochasticity enters via a low-dimensional and sufficiently simple probability distribution. In this paper we significantly extend the scope of the CSG method by presenting alternative ways to calculate the integration weights. A full convergence analysis for this new variant of the CSG method is presented and its efficiency is demonstrated in comparison to more classical stochastic gradient methods by means of a number of problem classes relevant to stochastic optimization and machine learning.
\end{abstract}

%\todoM{Abstract // Acknowledgement // Einleitung anpassen // Conclusion}
%\begin{itemize}
%    \item weiteres Setting als z.B. SG (Remark \ref{re:GeneralizedSetting})
%    \item Fehler im Gradienten geht gegen 0, d.h. mit dem Gradienten kann man auch was %anfangen, z.B. Abbruchkriterien
%    \item Approx. an Zielfunktionswert gibt's mit dazu
%    \item Gewichtsberechnungen (insbes. der inexakte hybrid) erlauben es, die Kosten pro Schritt quasi kontinuierlich gegen die Anzahl der Gesamtschritte zu tauschen. Teures Problem $\rightarrow$ lieber mehr Kosten in Gewichte und weniger Gesamtschritte, billiges Problem umgekehrt
%    \item höherdim. / andere Verteilungen / ...
%    \item Maß unbekannt / Daten
%    \item Anwendungen (mehr Breite)
%\end{itemize}
\section{Introduction}
%In this contribution we propose an efficient algorithm to solve deep learning tasks in the case of expensive (e.g. numerical complexity, data mining costs) generation of input-output data pairs.
%To solve problems in machine learning or deep learning, stochastic optimization schemes are commonly used. However, schemes based on the concept of \textit{stochastic gradient method} (SG) \cite{Monro1951} and its modification \textit{stochastic average gradient method} (SAG) \cite{LeRoux2017} do not rely on the Lipschitz continuity of the gradient. The presented algorithm can be seen as a generalization of these schemes. Depending on the probability measure and the choice of the underlying metric, the two mentioned schemes can be realized. However, the more interesting and promising case lies in between. This will be discussed in Section TODO on the basis of various examples.
%\todoM{Weiteres Setting stärker erwähnen / Verweis auf Anwendungsbeispiele weiter hinten / Andere Verfahren in Einleitung stärker beleuchten (A) / Machine Learning}
In the past, a variety of stochastic optimization schemes have been developed, e.g., \cite{TRMeth,LevelSetMeth,SARAH,InexProxMeth,PrimDualMeth}, in the context of optimization problems in which the expected-value of a cost function $j$ is minimized, i.e.,
\begin{equation}\label{eq:BspProblemEinleitung}
    \begin{aligned}
        \min_{\theta \in \Theta^{\textnormal{ad}}}\quad & \Ew[j(\theta,X)] = \int_{\mX}j(\theta,x)\mu(\mathrm{d}x),
    \end{aligned}
\end{equation}
 with probability measure $\mu$ and the associated random variables $\X$. Among the most popular algorithms are the \textit{stochastic gradient method} (SG) \cite{Monro1951} and its modification the \textit{stochastic average gradient method} (SAG) \cite{LeRoux2017}. Both of these methods have been analyzed extensively in the literature and are characterized by a low cost per iteration.

Nonetheless, SG and SAG suffer from a number of known disadvantages, like the lack of efficient stopping criteria (cf. \cite{Patel2016}) or optimal stepsize rules (cf. \cite{Paquette2020,tan2016barzilaiborwein}). To tackle these issues, a wide range of modified SG methods have been developed. For example, \cite{TRMeth} uses a trust-region-type model to normalize the steplengths, whereas the iSARAH algorithm proposed in \cite{SARAH} combines an inner SG scheme with an outer (inexact) full gradient descent method.

An additional disadvantage of SG is the relatively restrictive setting of \eqref{eq:BspProblemEinleitung}. In an effort to address this, \cite{InexProxMeth} and \cite{PrimDualMeth} suggest inexact proximal stochastic second-order methods and stochastic primal-dual fixed-point methods to allow for different types of objective functionals. For the case that the constraints include expected-valued functions, a level set method is also analyzed in \cite{LevelSetMeth}.

%Furthermore, especially in the context of large-scale optimization problems, where the evaluation of a gradient sample is time-consuming, the inherently ``wasteful'' nature of SG iterations can lead to poor performance.
The \textit{continuous stochastic gradient method} (CSG) proposed in \cite{pflug_CSG} is able to solve a wider class of problems. This is achieved by combining the information collected in previous iterations in an optimal way, meaning that CSG gains a significantly improved gradient approximation and is able to estimate the current objective function value during the optimization process. For a characterization of the class of problems that CSG can solve, we refer to \cref{re:GeneralizedSetting}. Here we simply note that these problems include objective functions with nested expectation values (\Cref{sec:ComparisonSCGD}) and problems with chance constraints (\Cref{sec:ChanceConstr}).

While these advantages of the original version of CSG \cite{pflug_CSG} are known, a serious drawback remains: to approximate function values and gradients as mentioned above, integration weights have to be computed via an analytical formula, which requires full knowledge of the probability measure $\mu$. Moreover, the evaluation is based on a Voronoi diagram, which is challenging to compute if the dimension of the parameter set $\mX$ is larger than 2. As a consequence, in \cite{pflug_CSG} only examples with a one-dimensional uniform distribution were presented.

In this contribution, we further expand the setting of CSG by introducing new methods to calculate the weights used for the gradient and cost function value approximations. This enables us to apply the CSG algorithm to problems of higher dimension, to arbitrary measure $\mu$ appearing in \eqref{eq:BspProblemEinleitung} and even to problems where $\mu$ might be unknown, e.g., in a data-driven context. 

Depending on the concrete setting, i.e., depending on the dimensions of $\theta,x$ and on how time-consuming the evaluation of a gradient sample is, the proposed methods allow us to continuously trade weight-computation time and speed of convergence (w.r.t. number of gradient sample evaluations).

In this article we present a full convergence analysis for this extended CSG method. In particular, we show that the errors in both the gradient approximation and the objective function value approximation vanish as the number of steps increases. As a consequence, these values can be utilized, for instance, to apply stopping criteria based on first order optimality conditions. Moreover, this allows the potential  combination of the CSG method with slightly adapted step length strategies as known from the world of deterministic optimization methods, a topic we leave open for future research.

The remaining structure of the paper is as follows. In \Cref{sec:Definitions}, we detail the mathematical structure of the problems we would like to solve using the CSG method. In \Cref{sec:TheCSGAlg}, we present the CSG method with generalized weight computation. \Cref{sec:konvergenzanalyse} is devoted to the convergence analysis and in \Cref{sec:numres} we compare the generalized CSG method to more traditional SG-type algorithms using three different classes of test problems.

%##################################################################################################
%################################### Section 2 ####################################################
%##################################################################################################

\section{Problem setting and definitions}\label{sec:Definitions}

Following the classic setup for expected-value objective functions, we introduce the set of admissible designs $\mP\subset\R^\dpar$ and the parameter set $\mX\subset\R^\din$, where $\dpar,\din\in\N$.
In the optimization process, the random samples $x_1,x_2,\ldots$ drawn from the parameter set $\mX$ are assumed to be realizations of independent uniformly random variables $X_i\sim\mu$ for all $i\in\N$, i.e., $X_1,X_2,\ldots$ are independent and follow an underlying probability distribution $\mu$, which may be unknown.

We define the following probability space setup:

\begin{definition}[Probability space setup] 
The probability space $(\Omega,\mathcal A,\Prob)$ is given by
\begin{equation*}
    \begin{aligned}
    \Omega &:= \mX^{\mathbb N}, \Prob := \mu^{\otimes \mathbb N},\\
    \mathcal A &:= \sigma(\{A_1\times \ldots\!\times\! A_n: A_i \in \mathcal B(\mX), \forall i, n \in \mathbb N\}),
    \end{aligned}
\end{equation*}
where  $\mu^{\otimes \mathbb N} (A_1\times \ldots \times A_n) = \prod_{i = 1}^n \mu(A_i)
$
is the product measure, $\mu$ is a probability measure on $\mX$ and $\sigma(\cdot)$ is the smallest $\sigma$-field that contains $\cdot$. We denote by $\supp(\mu)$ the support of the measure $\mu$, i.e.,
$$
\supp(\mu) := \{\x \in \mX:\ \mu(B_\eps(\x)) > 0 \ \forall \eps > 0\},
$$
where $B_\eps(x)$ denotes an open ball of radius $\eps > 0$ around $x \in \mathbb \mX$. We write $\X_n: \Omega \to \mX$, $(\omega_k)_{k \in \mathbb N} \mapsto \omega_n$ for the projection to the $n = 1,2,\ldots$ coordinate and define $\X:= \X_1$.
\end{definition}
With this setup, the objective function takes the following form:

\begin{definition}[Objective function] 
The objective function $J:\mP \rightarrow \R$ is given by
\begin{equation*}
J(\u) := \Ew\big[j(\theta,\X)\big] = \int_{\mX} j(\theta,\x) \, \mu(\mathrm d\x),
\end{equation*}
with a measurable function $j\in C^1(\mP \times \mX ; \R )$ and random variable $X$. 
\end{definition}
\begin{remark}[Generalization of the setting]\label{re:GeneralizedSetting}
During the optimization process, we may also generate an approximation $\hat{J}_n$ to the exact objective function value $J(\theta_n)$ with almost no additional computational cost. We will later show that $\Vert\hat{J}_n-\nabla J(\theta_n)\Vert_{\mP}\to 0$ (see \Cref{rem:ObjFuncApprox}).

This enables us to solve a much broader class of optimization problems where the objective function may depend non-linearly on the expression above, i.e.,
\begin{equation*}
    \tilde{J}(\theta) := f\big(\theta,\Ew[j(\theta,X)]\big),
\end{equation*}
with a Lipschitz continuously differentiable function $f:\mP\times\R\to\R$. Included in the set of possible objective functions are e.g. tracking functionals
\begin{equation*}
    \tilde{J}(\theta) := \tfrac{1}{2}\big\Vert h(\theta,\cdot) - f(\theta,\Ew[j(\theta,\cdot,X)])\big\Vert_{L^2}^2
\end{equation*}
and nested expected values
\begin{equation*}
    \tilde{J}(\theta) := \Ew_{\mathcal{Y}}\big[ f(Y,\Ew_{\mX}[j(\theta,X)])\big].
\end{equation*}
Note that such settings cannot be solved by SG algorithms.
\end{remark}
As we are aiming for a gradient based optimization scheme, we further state the derivative of the objective functional:

\begin{lemma}[Derivative of objective function]
The gradient of the objective functional $J$ is given by
   $\nabla J(\u) = \Ew\big[ \Delta(\u,\X)\big]$,
 where $\X \sim \mu$ and $\Delta: \mP\times \mX \rightarrow \R^\dpar$ denotes $\nabla_1j(\theta,\x)$.
\end{lemma}
\begin{proof}
This is a direct consequence of the linearity of the expectation value and the finite-dimensional derivative of $j$. Integration and differentiation can be exchanged due to the Lipschitz continuity of the integrand w.r.t. the integration variable.
\end{proof}

To state and prove convergence results for the algorithm presented in this work, we define the norms on the used spaces as follows:

\begin{definition}[Norms on $\mX,\mP$ and $\mP\times \mX$]\label{defi:norm}
In this contribution, we will use the notation $\|\cdot\|_\nX$ and $\|\cdot\|_\nT$ for the norm on the underlying spaces of the parameter space $\mX$ and the design space $\mP$, respectively. Due to norm-equivalence in finite dimensional spaces, the norm used in the spaces $\mP,\mX$ does not have to be specified and can be chosen according to the specific problem.
In addition, we define on $\mP\times \mX$ the following metric:
\begin{equation*}
d\big((\u,\x),(\hat \u ,\hat \x)\big) := \big\| \big( \|\u - \hat \u\|_\nT \, ,\, \|\x-\hat\x\|_\nX\big)\big\|_1 \ \forall (\u,\hat \u,\x,\hat \x) \in \mP^2 \times \mX^2.
\end{equation*}
 Choosing the 1-norm in the three-dimensional space as``outer"-norm is arbitrary and could for instance - in the other extreme case - be the $\infty$-norm. Of course this could also include positive weights for each individual component.
\end{definition}

\begin{assumption}[Regularity of the $\Delta$]\label{ass:reg_derivative}
 We assume $\Delta: \mP \times \mX \rightarrow \R^{\dpar}$ to be bounded and Lipschitz continuous, i.e., there exist constants $\CD,\LD\in \R_{>0}$ s.t. 
 \begin{equation*}
 \begin{aligned}
      \big\|\Delta(\theta,\x)\big\|_\nT &\leq \CD \\
     \big\|\Delta(\theta,\x) - \Delta(\hat\theta,\hat \x)\big\|_{\nT} &\leq \LD \big( \|\u-\hat\u\|_\nT+\|\x-\hat \x\|_\nX\big)
 \end{aligned}
 \end{equation*}
 for all $\u,\hat \u \in \mP$ and $\x,\hat\x \in \mX$. A sufficient condition for the stated regularity of $\Delta$ is thus to assume $\nabla j$ to be Lipschitz continuous in both arguments. 
\end{assumption}
 
The following  assumptions on the sets $\Uad, \Vad$ and the measure $\mu$ are an important ingredient for the convergence analysis of \cref{alg:2}.

\begin{assumption}[Regularity of $\Uad$, $\Vad$ and the measure $\mu$] \label{vor:V}
The set $\Uad \subset \R^{\dpar}$ is compact and convex. $\supp(\mu) \subset \Vad$ with $\Vad \subset \R^{\din}$ is open and bounded, and there exists $\Me,\Mz,\mathsf{M}_{_3}> 0$ s.t. $\forall\, \eps\in(0,\mathsf{M}_{_3})$ there exists $\mX_\eps\subset \mX$ satisfying $\mu(\mX_\eps) \geq 1- \Me\eps$ and
\begin{equation*}
\inf_{x \in \Vad_\eps} \mu\big(B_\eps(x) \big) \ge \Mz \eps^{\din},
\end{equation*}
where $B_\eps(x) \subset \mX$ is an open ball with radius $\eps$ centered in $x\in\mX$.
\end{assumption}

\begin{remark}[Examples for \Cref{vor:V}] In most cases, the choice $\mX_\eps = \mX$ is suitable, for example when $\mathcal{X}$ satisfies the uniform cone condition (cf. \cite[Definitioin 4.8]{UniformCone}). However, there are cases where the possibility of choosing $\mX_\eps \subset \mX$ as part of the conditions of \Cref{vor:V} allows us to consider even more general measures and sets.

Therefore, as an example, 
let $\mX := \{1/n : n\in \N\}$ and $\mu := \sum_{k=1}^\infty 2^{-k} \delta_{k^{-1}}(s)$. Then, for $M_n := [n^{-1},1] \cap \mX$, it holds that
\begin{equation*}
    \mu(M_n) = \sum_{k=1}^n 2^{-k} = 1-2^{-n}.
\end{equation*}
Thus, for $\eps = 2^{-n}$ and $c'=1$ we obtain $\mu(M_n) \geq 1-c'\eps$ and $\inf_{x\in M_n} \mu(B_\eps(x)) \geq 2^{-n} = \eps$.
Since $\inf_{x\in \mX} \mu(B_\eps(x))  = \mu(B_\eps(0)) = 2^{1-2^n}$, there exists no $c>0$ such that $2^{1-2^n} \geq c 2^{-n} \ \forall n \in \N$. %the latter would is equivalent to $2^{n+1-2^n} \geq c \ \forall n \in \N$ which can not hold.

For a uniform distribution and for all $0<p\le\infty$, the open unit ball w.r.t. the $p$-Norm 
\begin{equation*}
    \mathcal{X}^{p} := \big\{ x \in \mathbb{R}^{\din} : \Vert x\Vert_p < 1\big\}
\end{equation*}
also satisfy \Cref{vor:V}. While the case $1\le p\le\infty$ allows for $\mathcal{X}^p_\varepsilon = \mathcal{X}^p$, for $0<p<1$ we first have to obtain $\mathcal{X}^p_\varepsilon$ by trimming off the spikes of $\mX$.
\end{remark}

%#####################################################################################################
%############################# Section 3 #############################################################
%#####################################################################################################
 
\section{The algorithm}\label{sec:TheCSGAlg}

Before stating the algorithm, we first define the projection operator that ensures the sequence of generated designs $(\u_n)_{n\in\N}$ is in the set $\Uad$. 
\begin{definition}[Orthogonal projection] We define the orthogonal projection -- in the sense of $\|\cdot\|_\nT$ -- onto the set $\mP$ as follows:
\begin{equation*}%\label{eq:projektion}
    \Proj(\u) := \argmin_{\hat\u\in\Uad} \big\|\u - \hat\u\big\|_\nT.
\end{equation*}
Note that the existence and uniqueness of  $\Proj$  is guaranteed by the projection theorem (see e.g. %Theorem 1.4.1 and Remark 1.4.1 in
\cite{Aubin}) building on the convexity of $\mP$ as assumed in \Cref{vor:V}.
\end{definition}

\begin{lemma}[Properties of $\Proj$] \label{lem:Projection Properties}
    Let $\P\subset\R^{\dpar}$ satisfy \Cref{vor:V}. Then the following holds for all $x,y\in\R^{\dpar}$ and $z\in\P$:
    \begin{enumerate}
        \item[\textnormal{(a)}] $(\Proj(x)-x)^T(\Proj(x)-z)\le 0$,
        \item[\textnormal{(b)}] $(\Proj(y)-\Proj(x))^T(y-x)\ge\Vert\Proj(y)-\Proj(x)\Vert_{\nT}^2\ge 0,$
        \item[\textnormal{(c)}] $\Vert \Proj(y)-\Proj(x)\Vert_{\nT}\le \Vert y-x\Vert_{\nT}$.
    \end{enumerate}
\end{lemma}
\begin{proof}
    A proof of (a) can be found in \cite[Thm. 1.4.1 (ii)]{Aubin}. (b) and (c) correspond to (iii) and (ii) in \cite[Prop. 1.4.1]{Aubin}, respectively.
\end{proof}

Given $\theta_1$, $n=1$ and a sequence $\x_1,\x_2,\ldots$ of inputs assumed to be realizations of the independent random variables $\X_1,\X_2,\ldots$ introduced in \Cref{sec:Definitions}, the CSG method for the (possibly unknown) measure $\mu$ is given in \Cref{alg:2}. 
% ALGORITHMUS #####################################################################################
\begin{algorithm} 
\caption{CSG method}
\label{alg:2}
\begin{algorithmic}[1]
\WHILE{Termination condition not met}
\STATE{Sample objective function (optional):\\
    $\quad j_n := j(\theta_n,\x_n)$}
\STATE{Sample gradient:\\
    $\quad g_n := \nabla_\theta j(\theta_n,\x_n)$}
\STATE{Calculate weights $\alpha_k$}
\STATE{Calculate search direction:\\
    $\quad \hat G_{n} := \sum_{k=1}^n  \alpha_k g_k \vphantom{\Big(}$}
\STATE{Approximation to objective function value (optional):\\
    $\quad \hat J_{n} := \sum_{k=1}^n  \alpha_k j_k \vphantom{\Big(}$} %\Comment{Objective functional (optional)}
\STATE{Choose stepsize $\tau_n$}
\STATE{Gradient step:\\
    $\quad \theta_{n+1} := \Proj\big(\theta_n - \tau_n \hat G_{n}\big) \vphantom{\Big(}$}
\STATE{Update index:\\
    $\quad n \leftarrow n+1 \vphantom{\Big(}$}
\ENDWHILE
\end{algorithmic}
\end{algorithm}
%##################################################################################################
\FloatBarrier
\subsection{Calculating the weights}\label{sec:weights}
The quality of the weights $\alpha_k$ appearing in \Cref{alg:2} greatly impacts the accuracy of the gradient approximation $\hat{G}_n$ and therefore directly influences the overall performance of the CSG method. At the same time, a more optimal computation of the weights might be time-consuming. Since the trade-off between the time spent calculating the weights and the time gained by performing fewer gradient evaluations is heavily problem-specific, we propose four different methods for the weight-calculation in the $n$th step:
\subsubsection*{Exact} Following an exact nearest neighbor approximation for the integral
\begin{equation*}
    \nabla J(\theta_n) = \int_\mX \nabla_\theta j(\theta_n,x)\mu(\mathrm{d}x),
\end{equation*}
for each $k = 1,\ldots,n$, we define the set
\begin{equation*}
    M_k := \left\{ x\in\mX: d\big( (\theta_n,x),(\theta_k,x_k)\big) < d\big( (\theta_n,x),(\theta_j,x_j)\big)\text{ for all }j\in\{1,\ldots,n\}\setminus\{k\}\right\},
\end{equation*}
i.e., the set of points $x\in\mX$ such that $(\theta_n,x)$ is closer to $(\theta_k,x_k)$ than to any other previous evaluation point. Assuming that the measure $\mu$ is known, we then set ${\alpha_k := \mu(M_k)}$. This method has been thoroughly analyzed in \cite{pflug_CSG} and yields the best possible approximation to the exact gradient. However, it is computationally infeasible for problems of high dimensions.
\subsubsection*{Empirical} Utilizing the properties of the empirical measure $\mu_n$ (see \Cref{rem:emp_dist}), we may replace the exact weights mentioned above with the empirical weights
\begin{equation*}
    \alpha_k := \frac{1}{n}\sum_{i=1}^n 1_{M_k}(x_i) = \mu_n(M_k) \approx \mu(M_k),
\end{equation*}
where $1_{M_k}$ denotes the indicator function of the set $M_k$. Note that the computation of the empirical weights requires no knowledge of $\mu$ and is also feasible for high-dimensional problems. However,  many samples $x_i$ are needed to approximate the exact gradient with high accuracy.
\subsubsection*{Exact hybrid} Assuming that the dimension of $\mX$ is much smaller than the dimension of $\mP$, we might treat the designs and parameters separately. Instead of $M_k$, we now consider the sets
\begin{equation*}
    \widetilde{M}_i = \big\{ x\in\mX : \Vert x-x_i \Vert_\nX \le \Vert x-x_j\Vert_\nX \text{ for all }j=1,\ldots,n\big\}, \quad i = 1,\ldots,n.
\end{equation*}
The $\alpha_k$ are now calculated as a combination of the empirical and exact method
\begin{equation}\label{eq:HybridWeights}
    \alpha_k := \sum_{i=1}^n 1_{M_k}(x_i)\mu(\widetilde{M}_i).
\end{equation}
\subsubsection*{Inexact hybrid} As for the exact weights, the calculation of the exact hybrid weights requires knowledge of $\mu$. If $\mu$ is unknown, we may replace the factor $\mu(\widetilde{M}_i)$ in \eqref{eq:HybridWeights} with an empirical approximation. Since this only requires samples of $X$, which we assume are abundantly available, we can control the quality of this approximation through the number of samples we draw. The inexact hybrid weights are therefore calculated as follows:
\begin{equation*}
    \alpha_k := \frac{1}{\lfloor n^\beta\rfloor}\sum_{i=1}^n 1_{M_k}(x_{j_i})\sum_{m=1}^{\lfloor n^\beta\rfloor}1_{\widetilde{M}_{j_i}}(x_m),
\end{equation*}
where $\beta\ge 1$, $\lfloor n^\beta\rfloor$ is the total number of samples we have drawn until step $n$ and $x_{j_1},\ldots,x_{j_n}$ denotes the samples at which $\nabla_\theta j(\theta,x)$ has been evaluated.

\begin{figure}
    \centering
    \includegraphics[scale=0.85]{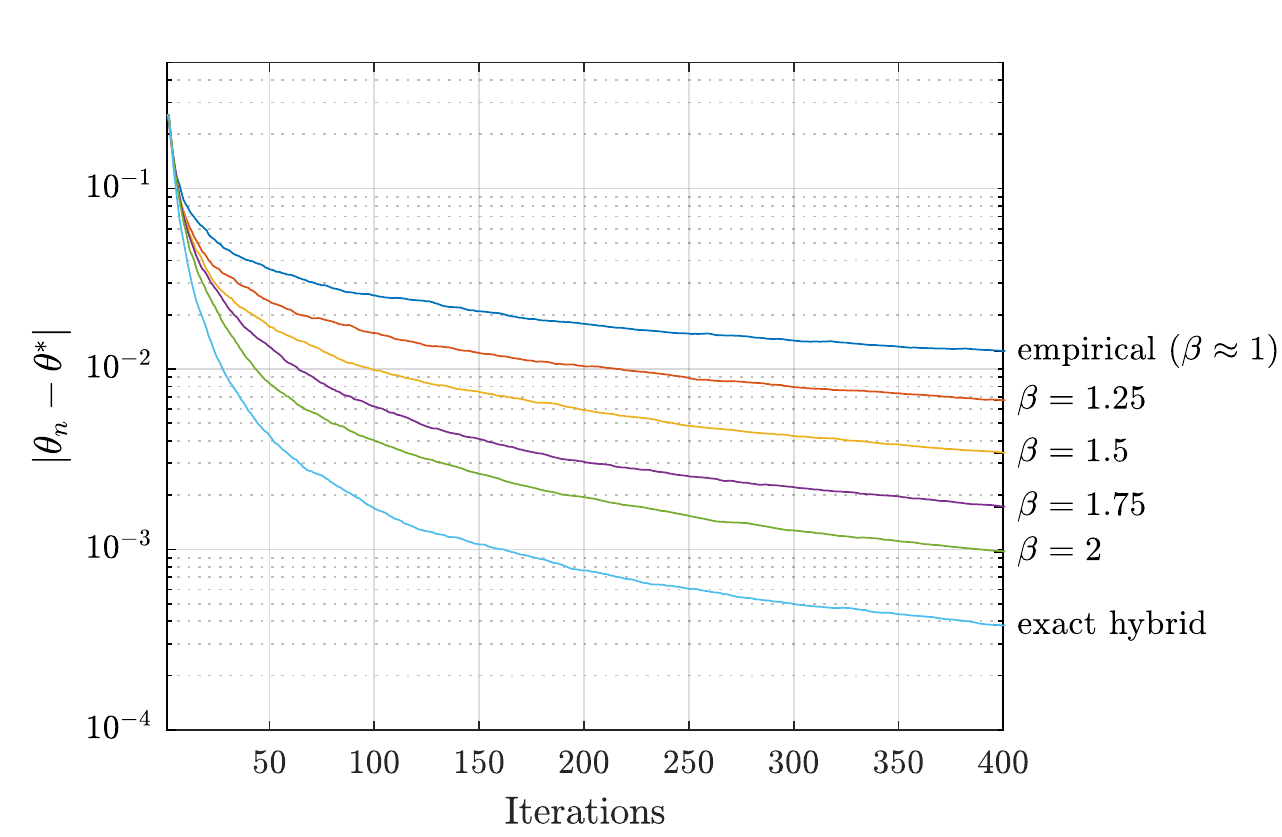}
    \caption{Absolute error $\vert \theta_n - \theta^\ast\vert$ in iteration $n$ for the setting $\mP = [-\frac{1}{2},\frac{1}{2}] = \mX$, $j(\theta,x) = \frac{1}{2}(\theta-x)^2$ and $X\sim\mathcal{U}_\mX$. The curves correspond to the median of 1000 runs with constant stepsizes $\tau_n = 1$ and randomized starting points in $\mP$.}
    \label{fig:Interpolation}
\end{figure}

\Cref{fig:Interpolation} shows that the inexact hybrid method allows us to interpolate between the purely empirical method and the exact hybrid variant by choosing an appropriate $\beta$.
\begin{remark}
In general, the nearest neighbor approximation used in all methods mentioned above worsens as the dimension of $\mP\times\mX$ increases (cf. \cite{NNMeaningful}). Especially for problems where $\dim(\mP)\ll\dim(\mX)$, Monte Carlo integration results (\cite{MonteCarloHighDim}) suggest that the performance boost gained by better weight calculation begins to become negligible. The proposed CSG methods are therefore best suited for optimization problems where $\mX$ is of small dimension in comparison to $\mP$ and the evaluation of $j(\theta,x)$ is time-consuming. 

Furthermore, to ensure the best possible performance, the metric $d$ should be chosen according to the specific problem.
\end{remark}

\begin{remark}[SAG and SG as two extreme cases of the algorithm]~\\
As stated in \Cref{defi:norm}, our metric $d$ can be chosen as
\begin{equation*}
    d\big((\theta,x),(\hat{\theta},\hat{x})\big) = a_1\Vert \theta-\hat{\theta}\Vert_\nT + a_2\Vert x - \hat{x}\Vert_\nX,
\end{equation*}
where $a_1,a_2 >0$ are arbitrary. By choosing $a_1 \gg a_2$, the nearest neighbor to $(\theta_n,x)$ is almost exclusively determined by the distance in the design variable. Hence, for the weights $\alpha_k$ we obtain $\alpha_n \approx 1$ and $\alpha_1,\ldots,\alpha_{n-1}\approx 0$, i.e., the CSG algorithm will behave very similarly to the usual SG algorithm.

Analogously, choosing $a_1 \ll a_2$ will lead to a performance similar to that of SAG.
\end{remark}

%##################################################################################################
%################################# Section 4 ######################################################
%##################################################################################################
\FloatBarrier
\section{Convergence analysis} \label{sec:konvergenzanalyse}
In this section we will study the convergence of the proposed algorithm. We will have to study probabilistic convergence behaviour in terms of ``almost sure convergence'', as we choose evaluation points randomly within the algorithm,. Therefore, we first state first order optimality conditions, assumptions on the regularity of the involved functions and on the steplength $\tau$, and a suitable probability space setting.
\subsection{Optimality conditions and assumptions}\label{subsec:1}
For $h \in C^1(\Uad)$ and $\Uad$ convex we have the following equivalent sufficient conditions for first order optimality:

\begin{corollary}[Optimality conditions]\label{kor:optimCond}
    For all $\u^* \in \mP$ the following items are equivalent:
    \begin{enumerate}
       \item[\textnormal{(a)}]$ -\nabla h(\u^*)^T(\u-\u^*) \leq 0 \quad \forall \u\in\Uad$
        \item [\textnormal{(b)}] 
        $\P(\u^* -t\nabla h(\u^*)) = \u^* \quad \forall t \geq 0$.
    \end{enumerate}
    A point $\u^*\in\Uad$ satisfying these conditions is called a stationary point.
\end{corollary}
\begin{proof}
  The proof can be found in e.g. \cite{pflug_CSG}. 
\end{proof}

To guarantee that \Cref{alg:2} generates a convergent subsequence, the stepsizes have to be damped, i.e., $(\tau_n)_{n \in \N}$ has to be a null series with upper and lower bound as stated in the following Assumption. However, in contrast to the ordinary stochastic gradient descent method, if \Cref{alg:2} generates -- with stepsizes satisfying $\tau_n \geq  \tau > 0\ \forall n\in \N$ -- a convergent sequence, the limit point is also a stationary point of the objective function. This is shown in \Cref{thm:Main theorem}.
\begin{assumption}[Steplength]\label{vor:schrittweitenfolge}
    The steplength $(\tau_n)_{n \in \N}$ in \Cref{alg:2} satisfies the following: $\exists N \in \N$, $\cu,\co\in \R_{>0}$ and $\D \in \big(0,\tfrac{1}{\max\{\din,2\}}\big)$ s.t. 
    \begin{equation*}
        \cu n^{-1}\leq \tau_n \leq \co { n^{-1+\frac{1}{\max\{\din,2\}}-\D}} \quad \forall n \in \N_{>N}.
    \end{equation*}
\end{assumption}
These bounds on the steplength satisfy the conditions stated in \cite[Eqns. (6) and (26)]{Monro1951}, as well as equivalently in \cite[Eqn. (4.19)]{Bottou2016} in the one-dimensional case, and can be seen as a higher dimensional equivalent.

In the following we assume that these assumptions are always satisfied without mentioning it explicitly.

\subsection{Error in the search direction}
In this subsection we analyse the error in the $n-$th iteration of the search direction $\hat G_n$ and the gradient of the objective functional $\nabla J_n$. To accomplish this, we define the following random variables:

\begin{definition}[Random variables]
For $\x\in\Vad$ and $\omega\in\Omega$ the sequence of random variables $\left(\Z_n\right)_{n\in \mathbb N}$ with $\Z_n: \Omega \times \mX \rightarrow \R_{\geq 0}$ is defined by
\begin{equation*}
    \Z_n(\omega,\x):= \min_{k = 1,\dots, n}d\big((\U_k(\omega),\X_k(\omega)),(\U_n(\omega),\x)\big),
\end{equation*}
where the designs $\U_k\in\Uad$ for $k>1$ depend by construction  on the initial design $\Theta_1$ and all ``previous''  random variables $\X_1(\omega),\ldots,\X_{k-1}(\omega)$, i.e., 
\begin{equation*}
    \Theta_k(\U_1, \X_1(\omega),\ldots,\X_{k-1}(\omega))
\end{equation*}
and is thus also a random variable. We shorten this dependency by the notation $\U_k(\omega)$.
\end{definition}
This random variable fulfills the following property:
\begin{lemma}\label{lem:konv1}
    For $\mu$ almost all $\x \in \supp(\mu)$ 
    \begin{equation*}
        \sum_{n = 1}^\infty \Prob\left(\Z_n(\cdot, \x) > \varepsilon_n \right) < \infty \quad \text{and} \quad 
        \sum_{n = 1}^\infty \sup_{\x \in \Vad_{\varepsilon_n} } \Prob\left(\Z_n(\cdot, \x) > \varepsilon_n \right) < \infty,
    \end{equation*}
    with \begin{equation}\label{eqn:eps_teps}
         \varepsilon_n := \tfrac{\CD\co}{1- 2^{-\frac{1}{2\max\{\din,2\}}}}\cdot n^{-\frac{\D}{2}} + \tilde \varepsilon_n
    \quad  \text{and} \quad
        \tilde \varepsilon_n := n^{\frac{\D}{2}-\frac{1}{2\max(\{\din,2\}}}.
    \end{equation} Here, $\CD$ is defined in \Cref{ass:reg_derivative}, $\Vad_{\varepsilon_n}$ in \Cref{vor:V} and $\co$,$\D$ in 
    \Cref{vor:schrittweitenfolge}.
\end{lemma}

\begin{proof}
We first define $i_0 \in \N$ as an auxiliary index as follows:
\begin{equation*}
    i_0 := \lceil n-a_n +1\rceil \quad \text{with} \quad a_n := 
     n^{1+\frac{\D}{2}-\frac{1}{\max\{\din,2\}}}. \label{eqn:def_an}
\end{equation*}
     By construction, we have
    \begin{equation*}
    \begin{aligned}
            \Prob(\Z_n(\cdot,\x) \geq \varepsilon_n) &\leq \Prob\left(\min_{k = i_0 ,\dots, n}d((\U_k,\X_k),(\U_n,\x)))  \geq \varepsilon_n \right) \\
            &\leq \Prob\left( \sum_{i=i_0}^{n-1}\|\tau_i \hat G_i  \|_\nT + \min_{k = i_0 ,\dots, n} \|\X_k - \x\|_\nX \geq \varepsilon_n \right) \\
            &\leq \Prob\left( \CD \sum_{i=i_0 }^{n-1} \tau_i + \min_{k = i_0 ,\dots, n} \|\X_k - \x\|_\nX  \geq \varepsilon_n \right).
    \end{aligned}
    \end{equation*}
   Observe that for $n > 2$ we obtain
    for all  $\kappa\in (0,1)$
    \begin{equation*}
    \begin{aligned}
        \sum_{i=i_0}^{n-1} \frac{1}{i^\kappa} &\le \int_{i_0-1}^n \tfrac{1}{s^\kappa} \, \mathrm d s = \frac{1}{1-\kappa}\cdot \left(n^{1-\kappa} - (\lceil n-a_n\rceil)^{1-\kappa} \right) 
        \le \frac{1}{1-\kappa}\cdot \left(n^{1-\kappa} - (n-a_n)^{1-\kappa} \right) \\
        &=  \frac{1}{1-\kappa}\cdot \left(\frac{n}{n^\kappa} - \frac{n-a_n}{(n-a_n)^{\kappa}} \right)
        = \frac{n}{1-\kappa}\cdot \left(\frac{(n-a_n)^\kappa - n^\kappa}{n^\kappa \cdot (n-a_n)^{\kappa}} \right) + \frac{a_n}{(1-\kappa)(n-a_n)^{\kappa}} \\
               &=  \frac{n}{1-\kappa}\cdot \left(\frac{n^\kappa (1-a_n/n)^\kappa - n^\kappa}{n^\kappa \cdot (n-a_n)^{\kappa}} \right) + \frac{a_n}{(1-\kappa)(n-a_n)^{\kappa}}.
    \end{aligned}
    \end{equation*}
Applying Bernoulli's inequality in the first term, we conclude
    \begin{align}
       \sum_{i=i_0}^{n-1} \frac{1}{i^\kappa} &\le \frac{n}{1-\kappa}\cdot \left(\frac{n^\kappa \left(1-\kappa\cdot \frac{a_n}{n}\right) - n^\kappa}{n^\kappa \cdot (n-a_n)^{\kappa}} \right) + \frac{a_n}{(1-\kappa)(n-a_n)^{\kappa}} \nonumber\\
        &= \frac{1}{1-\kappa}\cdot \left(\frac{-\kappa a_n}{(n-a_n)^{\kappa}} \right) + \frac{a_n}{(1-\kappa)(n-a_n)^{\kappa}}\nonumber\\
        &= \frac{a_n}{(n-a_n)^{\kappa}}
        = \frac{a_n}{n^\kappa}(1-\tfrac{a_n}{n})^{-\kappa}. \label{eq:BernoulliEq}  
    \end{align}
    Combining \Cref{vor:schrittweitenfolge} and \eqref{eq:BernoulliEq} yields
    \begin{equation*}
        \sum_{i=i_0}^{n-1} \tau_i \le \co\sum_{i=i_0}^{n-1} \frac{1}{i^{1+\D-\frac{1}{\max\{\din,2\}}}} \leq \co \tfrac{a_n}{n^\kappa}(1-\tfrac{a_n}{n})^{-\kappa},  
    \end{equation*}
    with $\D \in \left(0,\frac{1}{\max\{\din,2\}}\right)$ and  $\kappa := 1+\D - \frac{1}{\max\{\din,2\}} \in (0,1)$. 
    Hence, for $n\geq 2$ we obtain
    \begin{equation*}
        \left(1-\frac{a_n}{n}\right)^{-\kappa} = \left(1-n^{\frac{\D}{2}-\frac{1}{\max\{\din,2\}}}\right)^{-\kappa} 
        \leq \left(1-n^{-\frac{1}{2\max\{\din,2\}}}\right)^{-\kappa}
        \leq \left(1-2^{-\frac{1}{2\max\{\din,2\}}}\right)^{-1}.
    \end{equation*}
    Collecting these results, we see
    \begin{equation*}
        \sum_{i=i_0}^{n-1} \tau_i \le \frac{\co}{1- 2^{-\frac{1}{2\max\{\din,2\}}}} 
        \frac{a_n}{n^{1+\D - \frac{1}{\max\{\din,2\}} }} \leq \frac{\co}{1- 2^{-\frac{1}{2\max\{\din,2\}}}} n^{-\tfrac{\D}{2}}.
    \end{equation*}    
 Consequently,
 \begin{equation*}
     \sum_{i=i_0}^{n-1}\|\tau_i \hat G_i  \|_\nT \leq \frac{\CD\co}{1- 2^{-\frac{1}{2\max\{\din,2\}}}} n^{-\tfrac{\D}{2}} = \varepsilon_n-\tilde{\varepsilon}_n.
 \end{equation*}
    
    By \Cref{vor:V}, $\mu(\Vad\setminus \Vad_{\eps_n}) \rightarrow 0$. Hence, for $\mu$ almost all $\x \in \supp(\mu)$, there exists $n \in \N$ large enough, such that $\x \in \Vad_{\eps_n}$. Therefore,  
    \begin{equation*}
        \begin{aligned}
            &\Prob(\Z_n(\cdot,\x) \geq \varepsilon_n) \\
            &\quad\leq \Prob\left( \min_{k = i_0 ,\dots, n-1} \|\X_k - x\|_\nX \geq  \tilde \varepsilon_n \right) \\
            &\quad\leq \Prob\Big( \|\X_k - \x\|_\nX \geq  \tilde \varepsilon_n \ \forall k\in \left\{i_0 ,\dots, n-1\right\}\Big)\\
            &\quad=\prod_{k=i_0}^{n-1} \Prob\Big( \|\X_k - \x\|_\nX\geq  \tilde \varepsilon_n \Big)=\prod_{k=i_0}^{n-1} \left(1 - \mu\big(B_{\tilde \varepsilon_n }(\x) \big)  \right)\\
            &\quad \le\left(1 - \min\big\{\Mz (\tilde\varepsilon_n)^{\din} \, , \, 1 \big\} \right)^{\! a_n}.         
        \end{aligned}
    \end{equation*}
    As $\tilde \varepsilon_n \rightarrow 0$, there exists $N\in\N$ s.t. for $n\geq N$ we obtain
    \begin{equation*}
      \Prob(\Z_n(\cdot,\x) \geq \varepsilon_n)\leq\left(1 - \Mz n^{\frac{\din\, \D}{2}-\frac{\din}{2\max\{2,\din\}}}\right)^{\! a_n}.
    \end{equation*}

For simplicity, we define 
\begin{equation*}
   c_1 := \frac{\din\, \D}{2}-\frac{\din}{2\max\{2,\din\}}
\end{equation*} 
and recall that ${\log(1-x) \le -x}$ for all  $x < 1$. Since $c_1<0$, for $n$ large enough it holds
\begin{equation*}
\begin{aligned}
    &\left(1 - \Mz n^{\frac{\din\, \D}{2}-\frac{\din}{2\max\{2,\din\}}}\right)^{\! a_n}    
    = \left(1 - \Mz n^{c_1}\right)^{\! a_n} = \exp\left(a_n \log\left(1 - \Mz n^{c_1}\right)\right) \\
    &\le  \exp\left(-a_n \Mz n^{c_1}\right) 
    %&= \exp\left(- a_n  c_1 n^{c_2} \right) \\
    = \exp\left(- \Mz n^{1+\frac{ \D}{2}-\frac{1}{\max\{2,\din\}} + \frac{\din\, \D}{2}-\frac{\din}{2\max\{2,\din\}}} \right) \\
    &= \left\{\begin{array}{ll}
    \exp\left(- \Mz n^{1+\D-\frac{1}{2} -\frac{1}{4}} \right) & \din = 1 \\
    \exp\left(- \Mz n^{1+\frac{ \D}{2}-\frac{1}{\din} + \frac{\din\, \D}{2}-\frac{1}{2}} \right)& \din \ge 2 
    \end{array}\right. \ \le \exp\left(- \Mz n^{\D} \right). 
\end{aligned}
\end{equation*}
Recall that there is $N \in \mathbb N$ such that $\exp(-x) \le x^{-\frac{2}{\D}}$ for all $x \ge N$. It follows that for all sufficiently large $n$:
$\exp(-\Mz n^\D) \le \Mz^{-\frac{2}{\D}} n^{-2}$.
Hence,
\begin{equation*}
    \sum_{n = N}^\infty\left(1 - \Mz n^{\frac{\din\, \D}{2}-\frac{\din}{2\max\{2,\din\}}}\right)^{\! a_n}    \le \sum_{n = N}^\infty\exp\left(- \Mz n^{\D} \right) \le \sum_{n = N}^\infty  \Mz^{-\frac{2}{\D}} n^{-2} 
\end{equation*}
and thus 
\begin{equation*}
    \sum_{n = 1}^\infty \Prob\left(\Z_n(\cdot, \x) > \varepsilon_n \right) < \infty.
\end{equation*}
Finally, note that \Cref{vor:V} gives 
    \begin{equation*}
            \sup_{\x \in \Vad^{\varepsilon_n}} \Prob(\Z_n(\cdot,\x) \geq \varepsilon_n) \leq \left(1 -  c\cdot \tfrac{n^{\frac{\D}{2}}}{a_n}\right)^{\! a_n} = \left( 1-c\cdot n^{\frac{1}{\max\{\din,2\}}-1}\right)^{a_n} 
    \end{equation*}
with $c>0$. Following the same steps as above, we obtain
    \begin{equation*}
       \sum_{n = 1}^\infty \sup_{\x \in \Vad_{\varepsilon_n} } \Prob\left(\Z_n(\cdot, \x) > \varepsilon_n \right) < \infty. 
    \end{equation*}
\end{proof}
As a direct consequence of the latter result we obtain almost sure convergence. 
\begin{corollary}\label{cor:vanishingZ}
    For $\mu$ almost all $\x \in \supp(\mu)$
    \begin{equation*}
        \Z_n(\cdot, \x) \stackrel{a.s.}{\longrightarrow} 0 \quad\text{for}\quad n\to\infty.
    \end{equation*}
\end{corollary}
\begin{proof}
  The result follows by \Cref{lem:konv1} and the Borel-Cantelli Lemma (see for example Theorem 2.7 in \cite{klenke2013probability}). 
\end{proof}
\begin{remark}[Empirical distribution]\label{rem:emp_dist}
The empirical measure defined as 
\begin{equation}
\mu_n := \frac{1}{n}\sum_{i = 1}^n \delta_{\X_i}\label{eq:empiricalMeasure}
\end{equation}
satisfies
$\mu_n \Rightarrow \mu$ as $n \rightarrow \infty $
almost surely, see \cite[Theorem 3]{EmpiricalDistrPaper}. Here $\Rightarrow $ denotes the weak convergence of measures according to the weak-$^\ast$ convergence in dual space theory, i.e., 
\begin{equation*}
    \mu_n \Rightarrow \mu\quad \text{iff} \quad \int_\mX f(x) \mu_n(dx) \rightarrow \int_\mX f(x) \mu(dx) \quad \forall f \in C_b(\mathcal X,\mathbb R).
\end{equation*}

See for instance \cite{billingsley2013convergence} for the empirical distribution and  \cite[Section 7.3]{folland2009guide} for a functional analytical perspective on weak-* convergence in the function space setting discussed.

Since this property of $\mu_n$ is all we need in the following proofs and since the measures
\begin{equation*}
    \mu^{\textnormal{eh}}_n := \sum_{i=1}^n \delta_{X_i}\mu(\widetilde{M}_i)\quad \text{and}\quad \mu^{\textnormal{ih}}_n := \sum_{i=1}^n \delta_{X_{j_i}}\mu_{\lfloor n^\beta \rfloor}(\widetilde{M}_{j_i})\ ,
\end{equation*}
corresponding to exact hybrid weights and inexact hybrid weights respectively, also satisfy ${\mu^{\textnormal{eh}}_n \Rightarrow \mu}$ and ${\mu^{\textnormal{ih}}_n \Rightarrow \mu}$, we will w.l.o.g. work with empirical weights only. 
\end{remark}

Due to the Lipschitz continuity of $\Delta$ as defined in \Cref{ass:reg_derivative}, the expected value  $\nabla J(\u) = \Ew[\Delta(\u,\X)]$ is for $n \rightarrow \infty$ increasingly better approximated by  $\hat G_n$:

\begin{theorem}[Error in gradient approximation] \label{lem:approxkonvergenz}

    The norm of the difference between the search direction $\hat G_n$ and the gradient of the objective functional $\Ew[\Delta(\U_n,\X)]$ vanishes for $n\rightarrow \infty$, i.e.,  
    \begin{equation*}
        \| \hat{G}_n - \Ew[\Delta(\U_n,\X)]\|_\nT \; \stackrel{a.s.}{\longrightarrow} 0\quad\text{and} \quad \lim_{n\to\infty} \Ew\left[\|\hat{G}_n - \Ew[\Delta(\U_n,\X)]\|_\nT \right] = 0.
    \end{equation*}
\end{theorem}

\begin{proof}
   For $\x \in \supp(\mu)$ define 
   \begin{equation*}
  k^n(\omega;\x) := \argmin_{k = 1,\ldots,n} d\big((\U_k(\omega),\X_k(\omega)),(\U_n(\omega),\x)\big).
  \end{equation*}

For $\hat G_n$ as generated by \Cref{alg:2} with $n\in\N$ arbitrary but fixed, the following holds: 
     \begin{align}
      &\|\hat G_n - \Ew[\Delta(\U_n,\X)]\|_\nY \nonumber\\
      &\qquad= \bigg\|\ \sum_{i = 1}^n \int_{\mX} \delta_{k^n(\omega;\x)}(i)   \Delta(\U_i(\omega),\x_i) \mu_n(\mathrm d\x) - \int_{\mX}  \Delta(\U_n(\omega),\x) \mu(\mathrm d\x) \bigg\|_\nY\nonumber\\
      &\qquad\leq \bigg\|  \int_{\mX} \sum_{i = 1}^n \delta_{k^n(\omega;\x)}(i)   \Delta(\U_i(\omega),\x_i) - \Delta(\U_n(\omega),\x) \mu_n(\mathrm d\x) \bigg\|_\nY\nonumber\\
      &\qquad \qquad + \bigg\|\int_{\mX}  \Delta(\U_n(\omega),\x) \mu_n(\mathrm d\x) - \int_{\mX}  \Delta(\U_n(\omega),\x) \mu(\mathrm d\x) \bigg\|_\nY\nonumber\\
&\qquad\leq \LD  \int_{\mX} \Z_n(\omega,\x)\mu_n(\mathrm d\x)+ \bigg\|\int_{\mX}  \Delta(\U_n(\omega),\x) \mu_n(\mathrm d\x) - \int_{\mX}  \Delta(\U_n(\omega),\x) \mu(\mathrm d\x) \bigg\|_\nY, \label{eq:ZerlegungVonGradientenDifferenz}
  \end{align}
  where $\mu_n$ is the empirical measure given in \eqref{eq:empiricalMeasure} and $\LD$ the Lipschitz constant defined in \Cref{ass:reg_derivative}. We need to prove that both terms in \eqref{eq:ZerlegungVonGradientenDifferenz} vanish for $n \rightarrow \infty$. 
  
For the first term, the uniform (in $n$) Lipschitz continuity of $Z_n(\omega,\cdot)$ yields
\begin{equation*}
    \begin{aligned}
        \int_{\mX} Z_n(\omega,\x)\mu_n(\mathrm d\x) &= \int_{\mX} Z_n(\omega,x)\mu(\mathrm d\x) + \int_{\mX} Z_n(\omega,x)\mu_n(\mathrm dx) - \int_{\mX} Z_n(\omega,x)\mu(\mathrm d x) \\
            &\le \int_{\mX} Z_n(\omega,x)\mu(\mathrm d\x) + \LZ d_W(\mu_n,\mu),
    \end{aligned}
\end{equation*}
where $d_W$ denotes the Wasserstein distance of the measure $\mu_n$ and $\mu$ (see \cite{gibbs2002choosing}). Since $\mX$ is bounded, \cite[Theorem 6]{gibbs2002choosing} ensures that the Wasserstein distance metrizices the weak topology on the set of probability measures on $\mX$. Since $\mu_n \Rightarrow \mu$ almost surely, this gives 
$d_W(\mu_n,\mu) \rightarrow 0$ almost surely. Furthermore, by \Cref{ass:reg_derivative}, there exists $C>0$ s.t. $0 \le Z_n \le C$. Using \Cref{cor:vanishingZ}, we obtain $Z_n(\omega,x)\to0$ for almost all $\omega\in\Omega$. Therefore, Lebesgue's dominated convergence theorem yields
\begin{equation*}
    \int_{\mX} Z_n(\omega,x)\mu(\mathrm dx) \to 0\quad\text{for almost all }\omega\in\Omega.
\end{equation*}

To show that the second part of \eqref{eq:ZerlegungVonGradientenDifferenz} vanishes, observe that
\begin{equation*}
\begin{aligned}
 &\left\|\int_{\mX}  \Delta(\U_n(\omega),\x) \mu_n(\mathrm d\x) - \int_{\mX}  \Delta(\U_n(\omega),\x') \mu(\mathrm d\x') \right\|_\nT \\
 &\quad = \left\|\int_{\mX\times \mX}  \Delta(\U_n(\omega),\x) - \Delta(\U_n(\omega),\x)  Q_n(d(\x,\x')) \right\|_\nT   \\
 &\quad \le \LD \int_{\mX\times\mX}\|\x-\x'\|_\nX Q_n(d(\x,\x')),
\end{aligned}
\end{equation*}

where $Q_n(\cdot \times\Vad) = \mu_n$ and $Q_n(\Vad\times \cdot) = \mu$ is an arbitrary but fixed coupling of $\mu_n$ and $\mu$. By taking the infimum of all such couplings, we again obtain the Wasserstein distance $d_W$ of the measure $\mu_n$ and $\mu$, i.e., 
\begin{equation}\label{eqn:wass_dist}
 \left\|\int_{\mX}  \Delta(\U_n(\omega),\x) \mu_n(\mathrm d\x) - \int_{\mX}  \Delta(\U_n(\omega),\x) \mu(\mathrm d\x) \right\|_\nT\le \LD d_W(\mu_n,\mu).
\end{equation}
By the same arguments detailed earlier, $d_W(\mu_n,\mu)\to 0$ almost surely. Combining all the above yields 
\begin{equation*}
      \|\hat G_n - \Ew[\Delta(\U_n,\X)]\|_\nT \rightarrow 0
  \end{equation*}
  almost surely. Since the above quantities are bounded, the almost sure convergence also implies the convergence in expectation via Lebesgue's dominated convergence theorem. 
\end{proof}

\begin{remark}\label{rem:ObjFuncApprox}
Due to the regularity of $J$, we can show 
\begin{equation*}
    \Vert \hat{J}_n - J(\theta_n)\Vert_{\mP} \to 0
\end{equation*}
analogously to the proof of \Cref{lem:approxkonvergenz}.
\end{remark}
\begin{theorem}[Sum of error in gradient approximation]
    \label{lem:approxkonvergenz1}
    The expectation value of the summed norm of the difference between the search direction $\hat G_n$ and the gradient of the reduced objective functional $\nabla J$ weighted by the respected stepsize $\tau_n$ vanishes for  $n\rightarrow \infty$, i.e.,  
    \begin{equation} \label{eq:reihen konvergenz}
        \sum_{n = 1}^{\infty} \tau_n \Ew\left[ \|\hat{G}_n - \Ew[\Delta(\U,\X)] \|_\nT\right] < \infty.
    \end{equation}
\end{theorem}

\begin{proof}
Recall from the proof of \Cref{lem:approxkonvergenz} that 
    \begin{align}
      \Ew\left[\|\hat G_n - \Ew[\Delta(\U_n,\X)]\|_\nT\right]&\leq \LD  \Ew\left[\int_{\mX} \Z_n(\omega,\x)\mu_n(\mathrm d\x)\right]\nonumber\\
      &\quad + \Ew\left[\left\|\int_{\mX}  \Delta(\U_n(\omega),\x) \mu_n(\mathrm d\x) - \int_{\mX}  \Delta(\U_n(\omega),\x) \mu(\mathrm d\x) \right\|_\nT\right]. \label{eq:ZerlegungSumOfGradientError}
  \end{align}
We start by deriving an upper bound for the first term on the right hand side of this inequality.
Recall the definition of $\tilde\eps_n$ in \Cref{lem:konv1}, i.e.,  
\begin{equation*}
    \tilde\varepsilon_n := n^{-\frac{1}{\max\{2,\din\}}+\frac{\D}{2}},
\end{equation*}  with $\D$ as defined in 
    \Cref{vor:schrittweitenfolge}. Then, analogue to \Cref{lem:konv1} (cf. the proof and the notation there), together with 
\begin{equation*} 
D:= \sup_{(\tilde \theta,\tilde x),(\hat \theta,\hat x)\in\mP \times \mX}d((\tilde \theta,\tilde x),(\hat \theta,\hat x)),
\end{equation*}
we obtain the following estimate:
\begin{equation*}
    \begin{aligned}
        \Ew\left[\int_{\mX}\Z _n(\cdot,\x)  \mu_n(\mathrm  d\x) \right]
        &= \Ew\left[\int_{\mX}\Z_n(\cdot,\x) 1_{Z_n(\cdot,x) \leq \tilde\eps_n}(x) + \Z _n(\cdot,\x) 1_{Z_n(\cdot,x) > \tilde\eps_n}(x)  \mu_n(\mathrm  d\x)\right] \\
         &\leq \tilde\eps_n + D\Ew\left[\int_{\mX} 1_{Z_n(\cdot,x) > \tilde\eps_n}(x)  \mu_n(\mathrm  d\x)\right]\\
        &\leq \tilde\eps_n + D\Ew\left[\int_{\mX} \prod_{k=1}^n 1_{d((\U_k(\cdot),\X_k(\cdot)),(\U_n(\cdot),x) ) > \tilde\eps_n} \mu_n(\mathrm  d\x) \right].
    \end{aligned}
\end{equation*}
Setting $i_0 := \lceil n-a_n +1\rceil$ as in the proof of \Cref{lem:konv1} yields
\begin{equation*}
        \Ew\left[\int_{\mX}\Z _n(\cdot,\x)  \mu_n(\mathrm  d\x) \right]  \leq \tilde\eps_n + D\Ew\left[\int_{\mX} \prod_{k=i_0}^n 1_{d((\U_k(\cdot),\X_k(\cdot)),(\U_n(\cdot),x) ) > \tilde\eps_n} \mu_n(\mathrm  d\x) \right].
\end{equation*}
Since $\mu_n$ is the empirical measure as defined in \eqref{eq:empiricalMeasure} and due to the linearity of $\Ew$, we obtain
\begin{equation*}
    \begin{aligned}
        \Ew\left[\int_{\mX}\Z _n(\cdot,\x)  \mu_n(\mathrm  d\x) \right]  &= \tilde\eps_n + \frac{D}{n}\sum_{i=1}^n \Ew\left[ \prod_{\substack{k=i_0 \\ k\neq i}}^n 1_{d((\U_k(\cdot),\X_k(\cdot)),(\U_n(\cdot),\X_i(\cdot)) ) > \tilde\eps_n}  \right] \\
           &= \tilde\eps_n + \frac{D}{n}\sum_{i=1}^n  \prod_{\substack{k=i_0 \\ k\neq i}}^n \Prob\big(d((\U_k(\cdot),\X_k(\cdot)),(\U_n(\cdot),\X_i(\cdot)) ) > \tilde\eps_n \big), 
    \end{aligned}
\end{equation*}
where we used the independency of all $(X_i)_{i\in\N}$.
Finally, applying Fubini's theorem results in
    \begin{equation}
       \Ew\left[\int_{\mX}\Z _n(\cdot,\x)  \mu_n(\mathrm  d\x) \right] = \tilde\eps_n + \frac{D}{n}\sum_{i=1}^n  \prod_{\substack{k=i_0 \\ k\neq i}}^n \int_\mX\Prob\big(d((\U_k(\cdot),\X_k(\cdot)),(\U_n(\cdot),x) ) > \tilde\eps_n\big)  \mu(dx). \label{eq:SummeProduktWahrsch}
   \end{equation}
Let $\mX_{\tilde\eps_n}\subset\mX$ be the set given in \Cref{vor:V}. Following the same argumentation as in the proof of \Cref{lem:konv1}, we obtain
\begin{equation*}
  \begin{aligned}
      &\int_{\mX}  \Prob\big(d((\U_k(\cdot),\X_k(\cdot)),(\U_n(\cdot),x) ) > \tilde\eps_n \big) \mu(d\x) \\
      &\quad = \int_{\mX\setminus \mX_{ \tilde\eps_n}}  \Prob\big(d((\U_k(\cdot),\X_k(\cdot)),(\U_n(\cdot),x) ) > \tilde\eps_n \big) \mu(d\x)\\
      &\qquad + \int_{\mX_{ \tilde\eps_n}}  \Prob\big(d((\U_k(\cdot),\X_k(\cdot)),(\U_n(\cdot),x) ) > \tilde\eps_n \big) \mu(d\x) \\
      &\quad \le c'  \tilde\eps_n + \sup_{\x \in \mX_{ \tilde\eps_n}}\Prob\big(d((\U_k(\cdot),\X_k(\cdot)),(\U_n(\cdot),x) ) > \tilde\eps_n \big)
      \le  c'  \tilde\eps_n + \left(1 - c\cdot \tfrac{n^{\frac{\D}{2}}}{a_n}\right)^{\! a_n},
  \end{aligned}
\end{equation*}
with $a_n$ defined as in \Cref{lem:konv1} and $c>0$.
Utilizing $\log(1-x)\leq -x$ for all $x<1$ shows that
\begin{equation*}
\begin{aligned}
    \left(1 -  c\cdot \tfrac{n^{\frac{\D}{2}}}{a_n}\right)^{\! a_n} &= \exp\left(a_n\log\big(1-c\cdot n^{\frac{1}{\max\{\din,2\}}-1}\big)\right)
        \le \exp\left( -ca_n\cdot n^{\frac{\D}{2}}\right) \\
        &= \exp\left( -c\cdot n^{1+\D-\frac{1}{\max\{\din,2\}}}\right)
        \le \exp\big( -c\cdot n^\D\big)
        \le \tilde\eps_n
\end{aligned}
\end{equation*}
for $n$ large enough. Therefore, we have
\begin{equation*}
    \int_{\mX}  \Prob\big(d((\U_k(\cdot),\X_k(\cdot)),(\U_n(\cdot),x) ) > \tilde\eps_n \big)\mu(\mathrm d x) \le (c'+1)\tilde\eps_n.
\end{equation*}
Inserting this into \eqref{eq:SummeProduktWahrsch} yields
  \begin{equation}
      \Ew\left[\int_{\mX}\Z _n(\cdot,\x)  \mu_n(\mathrm  d\x) \right] \leq \tilde\eps_n + D\big((c'+1)\tilde\eps_n\big)^{n-i_0+1}\le \bar{c}\cdot\tilde\eps_n \label{eq:ErsterSummandSchranke}
  \end{equation}
for $n$ sufficiently large and some $\bar{c}>0$.

Now, to bound the second term in \eqref{eq:ZerlegungSumOfGradientError}, recall from \eqref{eqn:wass_dist} that 
\begin{equation*}
 \Ew\left[\left\| \int_{\mX}\Delta(\u_n,\x) \mu_n(d\x)  -  \Delta(\u_n,\x) \mu(d\x) \right\|_\nT \right]\le \LD\cdot \Ew\left[d_W(\mu_n,\mu)\right],
\end{equation*}
where $d_W$ is the Wasserstein distance. By \cite[Thm. 1 for $q=3$ and $p=1$]{Fournier2015}, for all $\din \geq 1$ there exists $\tilde C(\din) \in \R_{>0}$
s.t.:
\begin{align}
 \Ew\left[d_W(\mu_n,\mu)\right] &\le \hat C(\din) \cdot M_3 \left\{\begin{array}{ll}
 n^{-\frac{1}{2}}, &\quad \din = 1,\\
 n^{-\frac{1}{2}} \log(1+n), &\quad \din = 2,\\
 n^{-\frac{1}{{\din}}}, &\quad \din \ge 3,\\
 \end{array}\right.\nonumber\\
& \leq \hat C(\din) \cdot M_3 \cdot n^{-\frac{1}{\max\{\din,2\}}}\log(1+n), \label{eq:ZweiterSummandSchranke}
\end{align} 
with 
$
    M_3 := \big(\int_\mX \|x\|^3_{\mX} \mu(dx)\big)^{_{\nicefrac{1}{3}}}
$.

Substituting \eqref{eq:ErsterSummandSchranke} and \eqref{eq:ZweiterSummandSchranke} into \eqref{eq:ZerlegungSumOfGradientError} yields
\begin{equation}
    \sum_{n = N}^{\infty} \tau_n \Ew\left[ \|\hat{G}_n - \Ew[\Delta(\u,\X)] \|_\nT\right]
        \le \bar{c}\sum_{n=N}^{\infty} \tau_n\tilde\eps_n + \hat{C}(\din)M_3\sum_{n=N}^{\infty}\tau_nn^{-\frac{1}{\max\{\din,2\}}}\log(1+n) \label{eq:SummeGradientenfehlerZusammengesetzt}
\end{equation}
for $N\in\N$ large enough. By \Cref{vor:schrittweitenfolge}, we have 
\begin{equation*}
    \tau_n\le\co n^{-1-\D+\frac{1}{\max\{\din,2\}}}.
\end{equation*}
When this is inserted into \eqref{eq:SummeGradientenfehlerZusammengesetzt}, we have
\begin{equation*}
   \sum_{n = N}^{\infty} \tau_n \Ew\left[ \|\hat{G}_n - \Ew[\Delta(\u,\X)] \|_\nT\right]  \le  \bar{c}\co\sum_{n=N}^{\infty} n^{-1-\frac{\D}{2}} + \hat{C}(\din)M_3\co\sum_{n=N}^{\infty}n^{-1-\D}\log(1+n),
\end{equation*}
showing that
\begin{equation*}
    \sum_{n = 1}^{\infty} \tau_n \Ew\left[ \|\hat{G}_n - \Ew[\Delta(\u,\X)] \|_\nT\right]  \le \infty.
\end{equation*}
\end{proof}
Before we can present our main result, we now collect a few auxiliary results.
\begin{lemma}[Collection of auxiliary results] \label{lem:SammlungAusAnderemPaper} \ \\
    \begin{enumerate}
    \item[\textnormal{(a)}] The objective functional value in iteration $n\in\N$ satisfies
        \begin{equation*}
            J_{n+1}-J_n\le -\tfrac{1}{\tau_n}\Vert \theta_{n+1}-\theta_n\Vert_{\nT}^2+\phi_n,
        \end{equation*}
        where $\phi_n:=\tau_n\Vert\nabla J_n-\hat{G}_n\Vert_{\nT}\Vert\hat{G}_n\Vert_{\nT}+\tau_n^2C\Vert\hat{G}_n\Vert_{\nT}^2$ and $C\ge0$ denotes a constant depending only on the Lipschitz constants and suprema of the functions involved.
    \item[\textnormal{(b)}] For $\phi_n$ as defined above, it holds that $\sum_{n=1}^\infty\Ew[\phi_n] <\infty.$
    \item[\textnormal{(c)}] For all $t\ge0$, we have
        \begin{equation*}
            \Vert\Proj(\theta_n-t\hat{G}_n)-\theta_n\Vert_{\nT}\le\tfrac{t}{\tau_n}\Vert\theta_{n+1}-\theta_n\Vert_{\nT}.
        \end{equation*}
    \end{enumerate}
\end{lemma}
\begin{proof}
Assertions (a), (b) and (c) correspond to Lemma 16, Corollary 17 and Lemma 18 in \cite{pflug_CSG}. Note that, by \Cref{lem:approxkonvergenz1}, the proofs provided in \cite{pflug_CSG} can also be carried over to our setting.
\end{proof}
\begin{theorem}[Main theorem]\label{thm:Main theorem}
Let $(\u_n)_{n\in \N}$ be generated by \Cref{alg:2} with weights calculated by one of the methods mentioned in \Cref{sec:weights}. Then there exists a sub-sequence $(\u_{n_k})_{k\in \N}$ converging to a stationary point, i.e.,
    \begin{equation*}
        \liminf_{n\to\infty} \Ew\left[\|\Proj(\u_n - t \nabla J_n) -\u_n\|_\nT^2\right] = 0 \quad \text{for all $t\ge 0$.}
    \end{equation*}
On the other hand, assume the time-step series $(\tau_n)_{n\in\N}$ satisfies 
$\tau_n \geq \tau$ for all $n \in \N$ and some $\tau > 0$.
Let further $(x_n)_{n\in\N}$ be dense in $\Vad$ and assume $(\theta_n)_{n\in \N}$ converges to $\theta^* \in \Uad$. Then $\theta^*$ is a stationary point of $J$, i.e. 
    \begin{equation*}
        \|\Proj(\theta^* - t \nabla J(\theta^*)) -\theta^*\|_{\nT}^2 = 0 \quad \text{for all $t\ge0$.}
    \end{equation*}   
\end{theorem}
\begin{proof}
To prove the first part, we show
    \begin{equation}
        \sum_{n=1}^\infty\tau_n\Ew\left[\Vert\Proj(\theta_n-t\nabla J_n)-\theta_n\Vert_{\nT}^2\right]<\infty\quad\text{for all $t\ge0$}. \label{eq:ErsterTeilMainTheorem}
    \end{equation}
By the assumed compactness of $\P$ and regularity of $J$, we have
    \begin{equation*}
        J_{\text{inf}}:=\inf_{\theta\in\P}J(\theta)>-\infty.
    \end{equation*} 
    For arbitrary $N\in\N$, \Cref{lem:SammlungAusAnderemPaper} (a) gives
    \begin{equation*}
        \begin{aligned}
            J_{\text{inf}}-J_1 \le \Ew[J_{N+1}]-J_1 
                = \sum_{n=1}^N\Ew[J_{n+1}-J_n]
                \le\sum_{n=1}^N\left(-\frac{1}{\tau_n}\Ew\Big[\Vert\theta_{n+1}-\theta_n\Vert_{\nT}^2\Big]+\Ew[\phi_n]\right).
        \end{aligned}
    \end{equation*}
    Rearranging terms and utilizing \Cref{lem:SammlungAusAnderemPaper} (b) yields
    \begin{equation}\label{eq:EstimateInConvergenceRe}
        \sum_{n=1}^\infty\frac{1}{\tau_n}\Ew\Big[\Vert\theta_{n+1}-\theta_n\Vert_{\nT}^2\Big]\le J_1-J_{\text{inf}}+\sum_{n=1}^\infty\Ew[\phi_n]<\infty.
    \end{equation}
    By \Cref{lem:Projection Properties} (c) and \Cref{lem:SammlungAusAnderemPaper} (c), we obtain
    \begin{equation*}
        \begin{aligned}
            \Vert\Proj(&\theta_n-t\nabla J_n)-\theta_n\Vert_{\nT}^2 \\
                &\le \Big(\Vert\Proj(\theta_n-t\hat{G}_n)-\theta_n\Vert_{\nT}+\Vert\Proj(\theta_n-t\nabla J_n)-\Proj(\theta_n-t\hat{G}_n)\Vert_{\nT}\Big)^2 \\
                &\le\left(\frac{t}{\tau_n}\Vert\theta_{n+1}-\theta_n\Vert_{\nT}+t\Vert\hat{G}_n-\nabla J_n\Vert_{\nT}\right)^2 \le \frac{2t^2}{\tau_n^2}\Vert\theta_{n+1}-\theta_n\Vert_{\nT}^2+2t^2\Vert\hat{G}_n-\nabla J_n\Vert_{\nT}^2,
        \end{aligned}
    \end{equation*}
    where we used Young's inequality in the last line.
    Therefore, it holds that
    \begin{equation*}
        \sum_{n=1}^\infty\tau_n\Ew\Big[\Vert\Proj(\theta_n-t\nabla J_n)-\theta_n\Vert_{\nT}^2\Big] \le 2t^2\sum_{n=1}^\infty\frac{1}{\tau_n}\Ew\Big[\Vert\theta_{n+1}-\theta_n\Vert_{\nT}^2\Big]+2t^2\sum_{n=1}^\infty\tau_n\Ew\Big[\Vert\hat{G}_n-\nabla J_n\Vert_{\nT}^2\Big].
    \end{equation*}
    \eqref{eq:ErsterTeilMainTheorem} now follows from \eqref{eq:EstimateInConvergenceRe} and \Cref{lem:approxkonvergenz1}. \\
For the second part, observe that convergence of $(\theta_n)_{n\in\N}$ and density of $(x_n)_{n\in\N}$ in $\Vad$ yield
    \begin{equation*}
        Z_n(x) \to 0 \quad\text{for all $x\in\Vad$}.
    \end{equation*}
Therefore, by steps similar to those performed in the proof of \Cref{lem:approxkonvergenz}, it holds that
\begin{equation*}
    \Vert \hat{G}_n - \nabla J_n\Vert_{\nT} \to 0,
\end{equation*}
where $\nabla J_n$ denotes $\nabla J(\theta_n)$.
Hence, for all $t\ge0$ we obtain
\begin{equation*}
    \begin{aligned}
    \|\Proj(&\theta^* - t \nabla J(\theta^*)) -\theta^*\|_{\nT}^2 \\
        &= \lim_{n\to\infty}\Vert\Proj(\theta_n-t\nabla J(\theta_n))-\theta_n\Vert_{\nT} \\
        &\le \lim_{n\to\infty}\left( \Vert\Proj(\theta_n-t \hat{G}_n)-\theta_n\Vert_{\nT} + \Vert\Proj(\theta_n-t \hat{G}_n)-\Proj(\theta_n-t\nabla J_n)\Vert_{\nT}\right) \\
        &\le \lim_{n\to\infty}\left( \frac{t}{\tau_n}\Vert \theta_{n+1}-\theta_n\Vert_{\nT} + \Vert\Proj(\theta_n-t \hat{G}_n)-\Proj(\theta_n-t\nabla J_n)\Vert_{\nT}\right) \\
        &\le \lim_{n\to\infty} \frac{t}{\tau}\Vert \theta_{n+1}-\theta_n\Vert_{\nT} + \lim_{n\to\infty}t\Vert\hat{G}_n-\nabla J_n\Vert_{\nT} 
        = 0,
    \end{aligned}
\end{equation*}
where \Cref{lem:SammlungAusAnderemPaper} (c) was used for the second inequality.
\end{proof}

%##############################################################################################
%####################################### Section 5 ############################################
%##############################################################################################

\section{Numerical Results}\label{sec:numres}
In this section, we compare the CSG methods to suitable algorithms from the literature in three different settings. The comparison is based on the number of gradient evaluations, since these represent the time-consuming computations in complex optimization tasks.
\FloatBarrier
\subsection{Comparison with SG}
To start our numerical analysis, we consider the problem
\begin{equation*}
    \min_{\theta\in\mP}\quad \frac{1}{2}\int_\mX (x-\theta)^2\mathrm{d}x,
\end{equation*}
where $\mP = \mX = [-\frac{1}{2},\frac{1}{2}]$.

To study the behavior of our algorithm, we choose four different stepsizes ($n^{-1}$, $n^{-2/3}$, $n^{-1/3}$ and a constant stepsize of $1$) and track the absolute error in each iteration $\vert \theta_n - \theta^\ast\vert$. To ensure meaningful results, the 10000 starting points were chosen randomly in $\mP$. For comparison, we do the same for the ordinary stochastic gradient descent method (SG), since it is commonly used for problems resembling our example.
\begin{figure}[h]
    \centering
    \includegraphics[width=\textwidth]{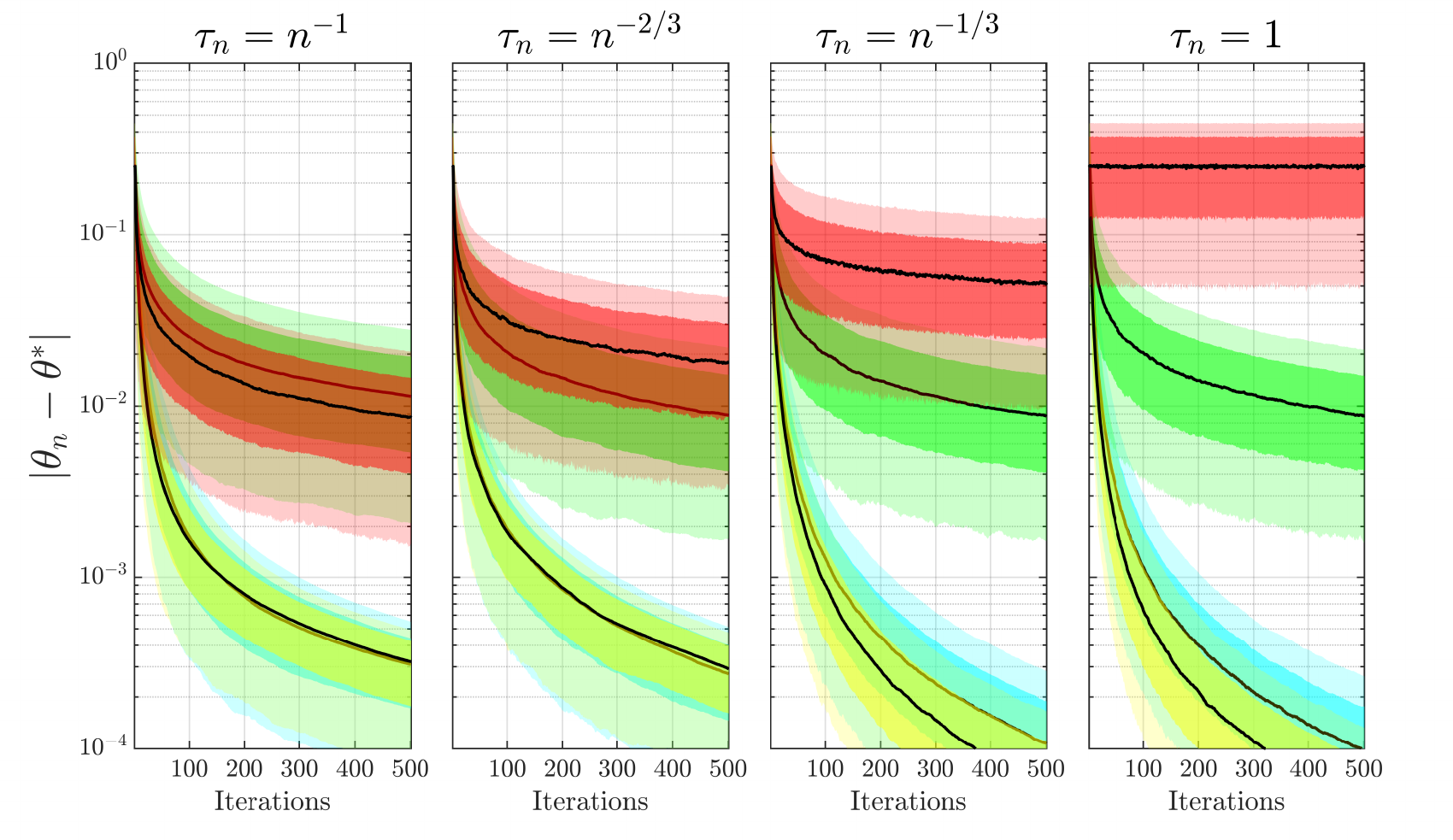}
    \caption{Comparison of the absolute error $\vert \theta_n-\theta^\ast\vert$ for SG (red), CSG with empirical weights (green), exact hybrid CSG (cyan) and exact CSG (yellow).}
    \label{fig:Beispiel1}
\end{figure}

Note that, in contrast to SG, a larger stepsize does not worsen the performance of the CSG algorithms for our example. A constant stepsize leads to a faster convergence for the hybrid and exact CSG method, whereas SG fails to solve the problem.

\subsection{Comparison with SCGD}\label{sec:ComparisonSCGD}
As mentioned in \Cref{re:GeneralizedSetting}, the vanishing error in inner function value approximations allows us to solve optimization problems in which the cost function depends non-linearly on a suitable expectation value.
For instance, we may solve the problem
\begin{equation}
    \min_{\theta\in\mP}\quad \frac{1}{20}\int_{\mathcal{Y}}\left(2y+5\int_{\mX}\cos\Big(\frac{\theta-x}{\pi}\Big)\mathrm{d}x\right)^2\mathrm{d}y, \label{eq:ProblemVergleichSCGD}
\end{equation}
where $\mP = [0,10]$, $\mX = [-1,1]$ and $\mathcal{Y} = [-3,3]$. The optimal solution $\theta^\ast = \frac{\pi^2}{2}$ to this example can be found analytically. Setting 
\begin{equation*}
    f_y(t) := \frac{3}{10}(2y+t)^2 \quad \text{and}\quad g_x(t) := 10\cos\Big(\frac{\theta-x}{\pi}\Big),
\end{equation*}
problem \eqref{eq:ProblemVergleichSCGD} can be reformulated as 
\begin{equation}
    \min_{\theta\in\mP}\quad \mathbb{E}_{\mathcal{Y}}\big[ f_y(\mathbb{E}_\mX[g_x(\theta)])\big].\label{eq:ProblemVergleichSCGDUmformuliert}
\end{equation}
Since $f_y$ is non-linear, the SG algorithm cannot be used to solve \eqref{eq:ProblemVergleichSCGD}. Therefore, we compare our results with the so called stochastic compositional gradient descent (SCGD) method (see \cite{SCGDPaper}), which is specifically designed for problems of the form \eqref{eq:ProblemVergleichSCGDUmformuliert}. 

\begin{figure}
    \centering
    \begin{minipage}{0.48\textwidth}
        \centering
        \includegraphics[width = \textwidth]{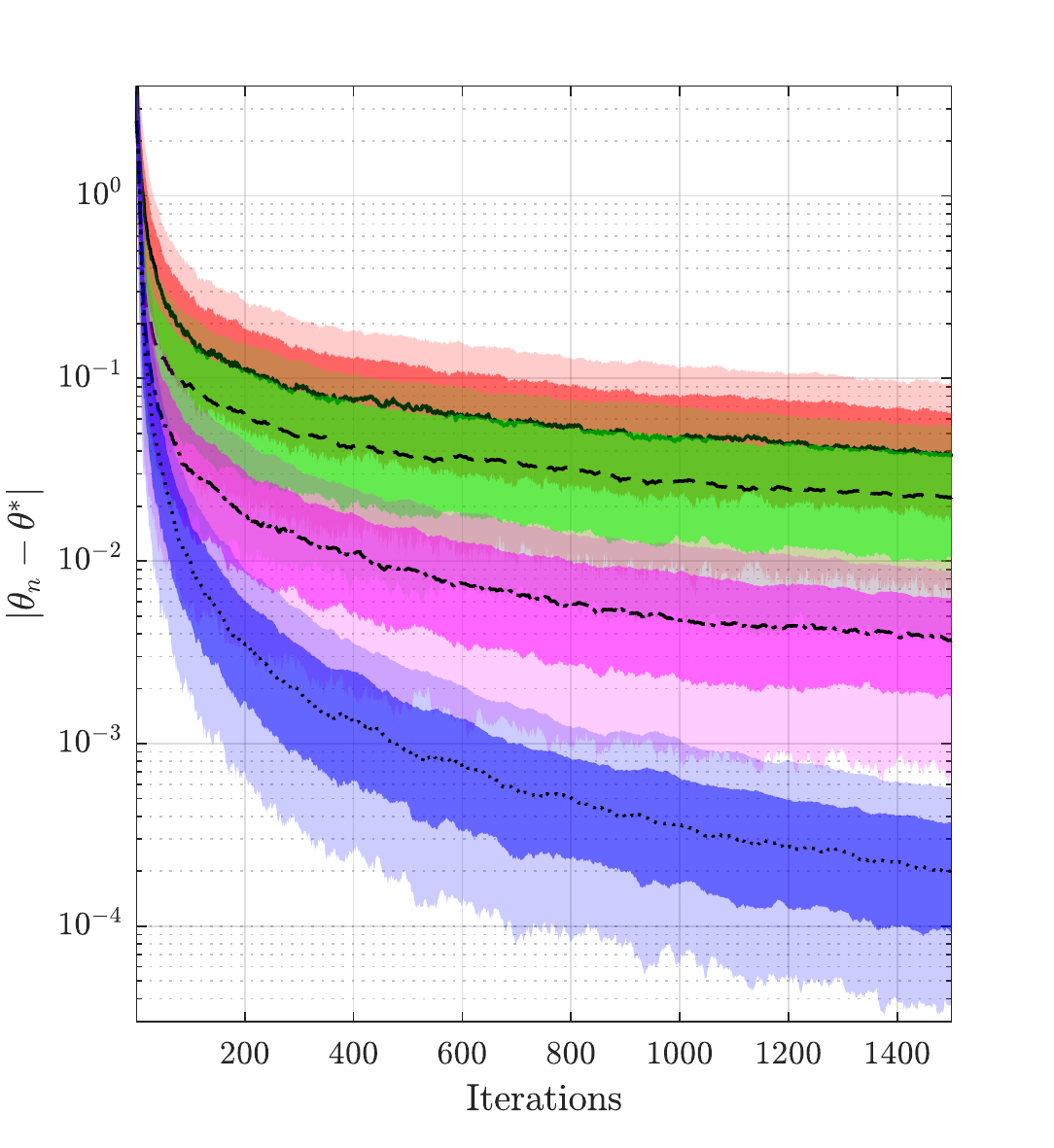}%
        \caption{Comparison of the absolute error $\vert\theta_n - \theta^\ast\vert$. From top to bottom: aSCGD (red/solid), CSG with empirical weights (green/dashed), inexact hybrid CSG with $\beta = 1.5$ (magenta/dash-dotted) and hybrid CSG (blue/dotted). The shaded areas indicate the quantiles $P_{0.1,0.9}$ (light) and $P_{0.25,0.75}$ (dark).}
        \label{fig:VergleichSCGD}
    \end{minipage}\hfill
    \begin{minipage}{0.48\textwidth}
        \centering
        \includegraphics[width = \textwidth]{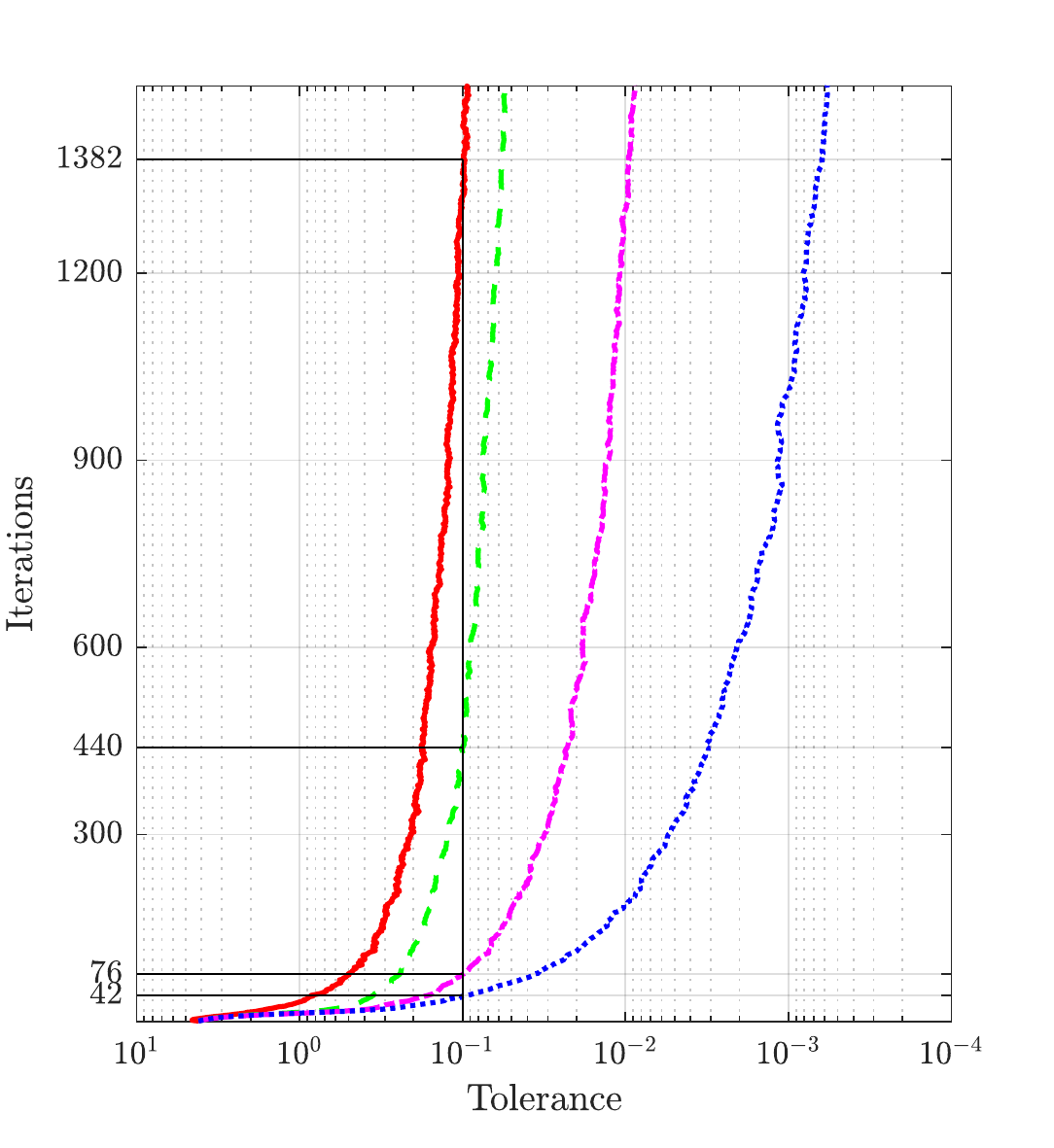}%
        \caption{Minimum number of steps needed for aSCGD (red/solid), CSG with empirical weights (green/dashed), inexact hybrid CSG with ${\beta = 1.5}$ (magenta/dash-dotted) and hybrid CSG (blue/dotted) such that at least 90\% of the runs achieved an absolute error smaller than the given tolerance.}
        \label{fig:StepsSCGD}
    \end{minipage}
\end{figure}

Again, the 1000 starting points are randomly generated. This time, however, we draw the starting points only from the interval $[\frac{11}{2},\frac{19}{2}]$ rather than $\mP=[0,10]$. The reason for this is that the optimal solution $\frac{\pi^2}{2}\approx 4.935$ would otherwise be very close to the median starting point, resulting in artificially small absolute errors for all methods.
Since the objective function in \eqref{eq:ProblemVergleichSCGD} is strongly convex in a neighborhood of the optimal solution, the accelerated SCGD method (see \cite{SCGDPaper}) performed better than the standard version. Therefore, we compared our results to the aSCGD algorithm and chose the optimal stepsizes for aSCGD according to Theorem 7 in \cite{SCGDPaper}. For the hybrid, inexact hybrid and empirical CSG algorithm, we chose a constant stepsize of $\frac{1}{30}$, which is a rough approximation to the inverse of the Lipschitz constant $\mathsf{L}_{\nabla J}$. The resulting graphs are shown in \Cref{fig:VergleichSCGD}.

From a practical viewpoint, the most interesting point is in how many iterations it takes the error to fall below a desired tolerance. For this purpose, we analyzed the number of steps after which the different methods achieved a given absolute error with 90\% certainty. The results can be seen in \Cref{fig:StepsSCGD}.

\subsection{Chance constraint problems}\label{sec:ChanceConstr}
\begin{figure}
    \centering
    \includegraphics[scale=0.85]{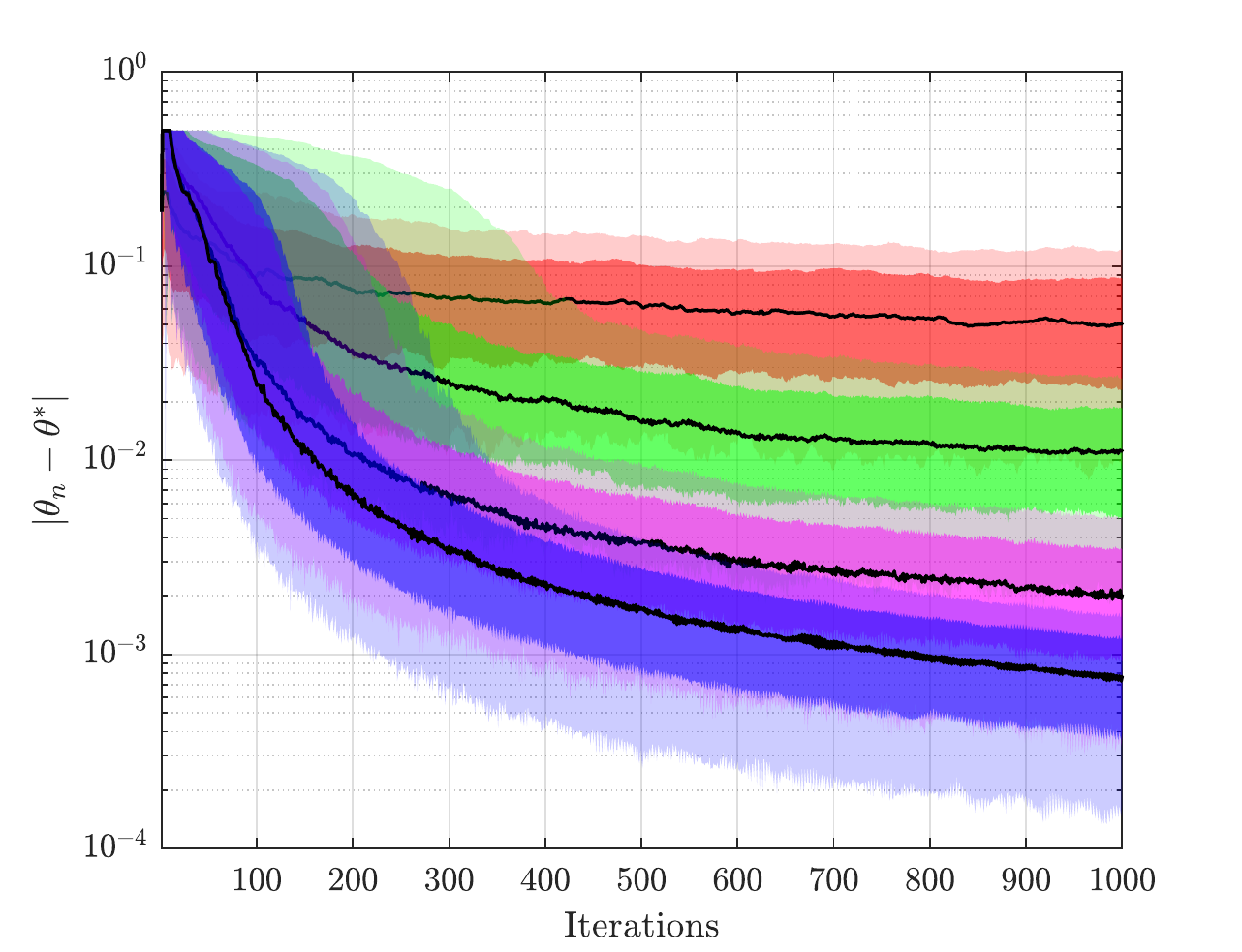}
    \caption{Comparison of the absolute error $\vert\theta_n-\theta_{opt}\vert$. From top to bottom: SCGD (red), CSG with empirical weights (green), inexact hybrid CSG with $\beta = 1.5$ (magenta) and hybrid CSG (blue). The shaded areas indicate the quantiles $P_{0.1,0.9}$ (light) and $P_{0.25,0.75}$ (dark).}
    \label{fig:ChanceConstraints}
\end{figure}

As a prototype example for chance constraint problems, we consider
\begin{equation*}
    \begin{aligned}
    \max_{\theta\in[0,\frac{3}{4}]}\quad &\theta \\
    \text{s.t.}\quad &\mathbb{P}(\theta-X^2\le 0) \ge \frac{1}{2}, \quad X\sim \mathcal{U}_{[-1,1]},
    \end{aligned}
\end{equation*}
with optimal solution $\theta^\ast=\frac{1}{4}$. By introducing the characteristic function $\chi_{[0,\infty)}$ and transforming the constraint to a penalty term, we arrive at
\begin{equation*}
    \max_{\theta\in[0,\frac{3}{4}]}\quad \theta - \lambda\max\left\{0,\frac{1}{2}\int_{-1}^1\chi_{[0,\infty)}(\theta-x^2)\mathrm{d}x - \frac{1}{2}\right\}.
\end{equation*}
Since the penalized objective function is no longer continuously differentiable, we cannot guarantee the existence of a gradient and will have to work with subgradients instead, cf. \cite{subgradients}. Note that the proofs provided above also hold true for a subgradient method if the stepsize is chosen accordingly.
While the computation of a (sub-)gradient of $\max\{0,\cdot\}$ is not an issue, $\chi_{[0,\infty)}$ needs to be regularized further. The final problem then reads as follows:
\begin{equation}
    \max_{\theta\in[0,\frac{3}{4}]}\quad \theta - \lambda\max\left\{0,\frac{1}{4}\int_{-1}^1\big((\tanh\big(\alpha(\theta-x^2)\big)+1\big)\mathrm{d}x - \frac{1}{2}\right\}. \label{eq:ChanceConstraintProblem}
\end{equation}
Due to the non-linearity of $\max\{0,\cdot\}$, we again choose the SCGD method for comparison. This time, the objective function is not strongly convex in a neighborhood of $\theta_{opt}$. Therefore, the stepsizes for the standard SCGD method are chosen according to Theorem 6 in \cite{SCGDPaper}, i.e., optimal for this setting. For the CSG algorithms, we choose $\tau_n = \frac{1}{n}$. Lastly, we fix $\lambda = 3$ and $\alpha = 25$. The optimal solution $\theta_{opt}$ to \eqref{eq:ChanceConstraintProblem} then satisfies $\vert\theta^\ast - \theta_{opt}\vert < 1.5\cdot 10^{-3}$. The results of 1000 runs with random starting points in $[0,\frac{3}{4}]$ are presented in \Cref{fig:ChanceConstraints}.

\FloatBarrier
\section{Conclusion and Outlook}
This article introduced a more flexible way to compute design dependent integration weights for the efficient approximation of the full cost function and its gradient when applying the CSG method to a class of stochastic optimization problems. While this approach significantly widened the scope of the CSG method, a number of interesting research questions remain to be investigated in the future. First, as a consequence of the strong convergence properties shown in this paper, the CSG method -- in the course of the optimization iterations -- behaves more and more like a fully deterministic descent method. This calls for more elaborate techniques to calculate the step length, e.g. linesearch or trust region strategies. Another interesting question is whether convergence of the iterates generated by the CSG method can be shown for a constant choice of step size. The numerical examples we have presented in this paper suggest that this should indeed be possible. And finally, exploiting specific structures of the given probability distributions could provide even more efficient integration techniques, enabling problems with high dimensional distributions to be solved more efficiently than using the empirical weight strategy presented in this article.

\section*{Acknowledgements}
Funded by the Deutsche Forschungsgemeinschaft (DFG, German Research Foundation) through project D05 in the CRC 1411 (Project-ID 416229255) and subproject B06 in TRR 154 (Project-ID 239904186).
\FloatBarrier
\bibliographystyle{siamplain}
\bibliography{references_siam}

\end{document}